\documentclass{article}[10pt]
\usepackage{a4}
\usepackage[T1]{fontenc}
\usepackage[latin1]{inputenc}
\usepackage[english]{babel}
\usepackage{hyperref}
\usepackage[dvips]{graphics}
\usepackage{amsmath}
\usepackage{amssymb}
\usepackage{amsopn}
\usepackage{amsbsy}
\usepackage{amstext}
\usepackage{latexsym}
\usepackage{amsfonts}
\usepackage{array}
\usepackage{epsfig}
\usepackage[mathscr]{eucal}
\title{Nontrivial lower bounds for the least common multiple of
some finite sequences of integers}
\author{Bakir FARHI}
\date{}

\def\ppcm{\mathrm{lcm}}
\def\pgcd{\mathrm{gcd}}
\let\epsilon=\varepsilon
\def\mod{\mathrm{mod}}

\let\m=\ppcm
\let\d=\pgcd

\newtheorem{thm}{Theorem}

\newtheorem{lemme}[thm]{Lemma}

\newtheorem{coll}[thm]{Corollary}
\newtheorem{conj}[thm]{Conjecture}

\begin{document}
\maketitle \vspace{-6cm}
\begin{flushleft}
{\it J. Number Theory}, \\
{\bf 125} (2007), p. 393-411.\\
N.B: This article gives the details of my earlier note \cite{f}.
\end{flushleft}~\vspace{3cm}

\begin{center}
Département de Mathématiques, Université du Maine, \\
Avenue Olivier Messiaen, 72085 Le Mans Cedex 9, France. \\
bakir.farhi@gmail.com
\end{center}
\vskip 1cm \noindent {\small {\bf Abstract.---} We present here a
method which allows to derive a nontrivial lower bounds for the
least common multiple of some finite sequences of integers. We
obtain efficient lower bounds (which in a way are optimal) for the
arithmetic progressions and lower bounds less efficient (but
nontrivial) for quadratic sequences whose general term has the
form $u_n = a n (n + t) + b$ with $(a , t , b) \in {\mathbb{Z}}^3
, a \geq 5 , t \geq 0 , \rm{gcd}(a , b) = 1$. From this, we deduce
for instance the lower bound: $\mathrm{lcm}\{1^2 + 1 , 2^2 + 1 ,
\dots , n^2 + 1\} \geq 0,32 (1,442)^n$ (for all $n \geq 1$).

In the last part of this article, we study the integer
$\mathrm{lcm}(n , n + 1 , \dots , n + k)$ $(k \in \mathbb{N} , n
\in {\mathbb{N}}^*)$. We show that it has a divisor $d_{n, k}$
simple in its dependence on $n$ and $k$, and a multiple $m_{n, k}$
also simple in its dependence on $n$. In addition, we prove that
both equalities: $\mathrm{lcm}(n , n + 1 , \dots , n + k) = d_{n,
k}$ and $\mathrm{lcm}(n , n + 1 , \dots , n + k) = m_{n, k}$ hold
for an infinitely many pairs $(n , k)$.} \vskip 2mm \noindent
\hfill{--------} 
\section{Introduction and notations}
In this article, $[x]$ denotes the integer part of a given real
number $x$. Further, we say that a real $x$ is a multiple of a
non-zero real $y$ if the ratio $x/y$ is an integer.

The prime numbers theorem (see e.g. \cite{har}) shows that
$\lim_{n \rightarrow + \infty} \frac{\log\mathrm{lcm}\{1, \dots,
n\}}{n } = 1$. This is equivalent to the following statement:
$$\forall ~\!\!\epsilon > 0 , \exists~\!\! N = N(\epsilon) ~\!/~\! \forall n \geq N :
~~~~ (e - \epsilon)^n ~\leq~ \ppcm\{1, \dots, n\} ~\leq~ (e +
\epsilon)^n .$$

Concerning the effective estimates of the numbers
$\mathrm{lcm}\{1, \dots, n\}$ $(n \geq 1)$, one has among others,
two main results. The first one is by Hanson \cite{han} which
shows (by using the development of the number $1$ in Sylvester
series) that $\mathrm{lcm}\{1 , \dots , n\} \leq 3^n$ for all $n
\geq 1$. The second one is by Nair \cite{nai} which proves (simply
by exploiting the integral $\int_{0}^{1 } x^n (1 - x)^n d x$) that
one has $\mathrm{lcm}\{1 , \dots , n\} \geq 2^n$ for all $n \geq
7$.

In this, we present a method which allows to find a nontrivial
lower bounds for the least common multiple of $n$ consecutive
terms $(n \in { \mathbb{N}}^*)$ of some sequences of integers. We
obtain efficient lower bounds (which in a way are optimal) for the
arithmetical progressions (see Theorem \ref{t4}). Besides, we also
obtain less efficient lower bounds (but nontrivial) for the
quadratic sequences whose general term has the form: $u_n = a n (n
+ t) + b$ with $(a, t, b) \in { \mathbb{Z}}^3$, $a \geq 5, t \geq
0, \pgcd(a, b) = 1$ (see Corollary \ref{c1}).

Our method is based on the use of some identities related to the
sequences which we study. More precisely, let ${(\alpha_i)}_{i \in
I}$ be a given finite sequence of nonzero integers. We seek an
identity of type $\sum_{i \in I} \frac{1}{\alpha_i \beta_i} =
\frac{1}{\gamma}$ where $\beta_i$ $(i \in I)$ and $\gamma$ are
nonzero integers. If $\mathrm{lcm}\{\beta_i, i \in I\}$ is bounded
(say by a real constant $R > 0$), one concludes that
$\mathrm{lcm}\{\alpha_i, i \in I \} \geq \frac{\gamma}{R}$ (see
Lemma \ref{l1}). It remains to check whether this later estimate
is nontrivial or not.

However, the point is that looking for identities of the above
types is not easy. Theorem \ref{t1} stems from concrete and
interesting example of such identities. Though, it is not likewise
that we can find other nontrivial applications, than the ones
presented here, for that specific example. In order to have
nontrivial lower bounds of least common multiple for other
families of finite sequences, it could be necessary to seek for
new identities related to those sequences.

In the last part of this article, we study the least common
multiple of some number of consecutive integers, larger than a
given positive integer. In Theorem \ref{t8}, we show that the
integer $\mathrm{lcm}\{n , n + 1 , \dots , n + k\}$ $(n \in
{\mathbb{N}}^* , k \in \mathbb{N})$ has a divisor $d_{n, k}$
simple in its dependence on $n$ and $k$ and a multiple $m_{n, k}$
simple in its dependence on $n$. In addition, we prove that $d_{n,
k}$ and $m_{n, k}$ are optimal in that sense that the equalities
$\mathrm{lcm}\{n, \dots, n + k\} = d_{n, k}$ and $\mathrm{lcm}\{n
, \dots , n + k\} = m_{n, k}$ hold for infinitely many pairs $(n,
k)$. More precisely, we show that both equalities are satisfied at
least when $(n, k)$ satisfies some congruence modulo $k!$ (see
Theorem \ref{t9}).
\section{Results}
\subsection{Basic Results}
\begin{lemme}\label{l1}
Let ${(\alpha_i)}_{i \in I}$ and ${(\beta_i)}_{i \in I}$ be two
finite sequences of non-zero integers such that:
$$\sum_{i \in I} \frac{1}{\alpha_i \beta_i} ~=~ \frac{1}{\gamma}$$
for some non-zero integer $\gamma$. Then, the integer
$\ppcm\{\alpha_i , i \in I\} . \ppcm\{\beta_i , i \in I\}$ is a
multiple of $\gamma$.
\end{lemme}
\begin{thm}\label{t1}
Let ${(u_k)}_{k \in \mathbb{N}}$ be a strictly increasing sequence
of non-zero integers. Then, for any positive integer $n$, the
integer:
$$\ppcm\left\{u_0 , \dots , u_n\right\} . \ppcm\left\{\prod_{
0 \leq i \leq n , i \neq j} \!\!\!(u_i - u_j) ~;~ j = 0 , \dots ,
n\right\}$$ is a multiple of the integer $(u_0 u_1 \dots u_n)$.
\end{thm}
\subsection{Results about the arithmetic progressions}
\begin{thm}\label{t2}
Let ${(u_k)}_{k \in \mathbb{N}}$ be a strictly increasing
arithmetic progression of non-zero integers. Then, for any
non-negative integer $n$, the integer $\ppcm\{u_0 , \dots , u_n\}$
is a multiple of the rational number:
$$\frac{u_0 \dots u_n}{n! \left(\pgcd\{u_0 , u_1\}\right)^n} .$$
\end{thm}
\begin{thm}[Optimality of Theorem \ref{t2}]\label{t3}
Let ${(u_k)}_{k \in \mathbb{N}}$ be a strictly increasing
arithmetic progression of non-zero integers such that $u_0$ and
$u_1$ are coprime. Then, for any positive integer $n$ which
satisfies:
$$u_0 u_n \equiv 0 ~\mod(n!) ,$$
we have:
$$\ppcm\{u_0 , \dots , u_n\} ~=~ \frac{u_0 \dots u_n}{n!} .$$
\end{thm}
\begin{thm}\label{t4}
Let ${(u_k)}_{k \in \mathbb{N}}$ be an arithmetic progression of
integers whose difference $r$ and first term $u_0$ are positive
and coprime. Then:
\begin{description}
\item[1)] For any $n \in \mathbb{N}$, we have:
$$\ppcm\{u_0 , \dots , u_n\} ~\geq~ u_0 (r + 1)^{n - 1} .$$
 Besides, if $n$ is a multiple of $(r + 1)$, we have:
$$\ppcm\{u_0 , \dots , u_n\} ~\geq~ u_0 (r + 1)^n .$$
\item[2)] For any $n \in \mathbb{N}$, we have:
$$\ppcm\{u_0 , \dots , u_n\} ~\geq~ r (r + 1)^{n - 1} .$$
\item[3)] For any $n \in \mathbb{N}$, we have:
$$\ppcm\{u_0 , \dots , u_n\} ~\geq~ \frac{n}{n + 1} r \left\{(r + 1)^{n - 1} + (r - 1)^{n - 1}\right\} .$$
\item[4)] For any $n \in \mathbb{N}$ satisfying $n \geq u_0 -
\frac{3 r + 1}{2}$, we have:
$$\ppcm\{u_0 , \dots , u_n\} ~\geq~ \frac{1}{\pi} \sqrt{r} (r + 1)^{n - 1 + \frac{u_0}{r}} .$$
\end{description}
\end{thm}
\begin{conj}[recently confirmed by S. Hong and W. Feng
\cite{hf}]\label{conj1}~\\ In the situation of Theorem \ref{t4},
we have for any $n \in \mathbb{N}$:
$$\ppcm\{u_0 , \dots , u_n\} ~\geq~ u_0 (r + 1)^n .$$
\end{conj}
The two following Theorems study the optimality of the part {\bf
4)} of Theorem \ref{t4}.
\begin{thm}\label{t5}
The coefficient $- \frac{3}{2}$ affected to $r$ which appears in
the condition ``$n \geq u_0 - \frac{3 r + 1}{2}$''$~\!\!$ of the
part {\bf 4)} of Theorem \ref{t4} is optimal.
\end{thm}
\begin{thm}\label{t6}~
\begin{description}
\item[1)] The optimal absolute constant $C$ for which the
assertion:
\begin{quote}
``For any arithmetic sequence ${(u_k)}_k$ as in Theorem
\ref{t4} and for any non-negative integer $n$ satisfying $n \geq u_0 - \frac{3 r + 1}{2}$, we have:\\
$\ppcm\{u_0 , \dots , u_n\} ~\geq~ C \sqrt{r} (r + 1)^{n - 1 +
\frac{u_0}{r}}$''
\end{quote}
is true, satisfies:
$$\frac{1}{\pi} ~\leq~ C ~\leq~ \frac{3}{2} .$$
\item[2)] More generally, given $n_0 \in \mathbb{N}$, the optimal
constant $C(n_0)$ (depending uniquely on $n_0$) for which the
assertion:
\begin{quote}
``For any arithmetic sequence ${(u_k)}_k$ as in Theorem \ref{t4}
and for any integer $n$ satisfying $n \geq \max\{n_0 , u_0 -
\frac{3 r + 1}{2}\}$, we have:
\\ $\ppcm\{u_0 , \dots , u_n\} ~\geq~ C(n_0) \sqrt{r} (r + 1)^{n -
1 + \frac{u_0}{r}}$''
\end{quote}
is true, satisfies:
$$\frac{1}{\pi} ~\leq~ C(n_0) ~<~ 4 (n_0 + 4) \sqrt{n_0 + 4} .$$
\end{description}
\end{thm}\newpage\noindent
{\bf Comments:}~\label{r1}
\begin{description}
\item[i)] The lower bound proposed by Conjecture \ref{conj1}
(recently confirmed in \cite{hf}) is optimal on the exponent $n$
of $(r + 1)$. Indeed, for any positive integer $n$ and for any
arithmetic progression ${(u_k)}_k$ as in Theorem \ref{t4}, we
obviously have:
$$\ppcm\{u_0 , \dots , u_n\} ~\leq~ u_0 u_1 \dots u_n ~\leq~ u_0 \left(\max\{u_0 , n\}\right)^n (r + 1)^n .$$
For any given positive real $\epsilon$, we can choose two
arbitrary positive integers $u_0$ and $n$ and a positive integer
$r$, which is coprime with $u_0$ and sufficiently large as to have
$(r + 1)^{\epsilon} > \left(\max\{u_0 , n\}\right)^n$. The
arithmetic progression ${(u_k)}_k$, with first term $u_0$ and
difference $r$, will then satisfy:
$$\ppcm\{u_0 , \dots , u_n\} ~<~ u_0 (r + 1)^{n + \epsilon} .$$
\item[ii)] A similar argument to that of the above part {\bf i)}
shows that the exponent $(n - 1)$ of $(r + 1)$ which appears in
the lower bound of the part {\bf 2)} of Theorem \ref{t4} is
optimal.
 \item[iii)] For small values of $n$ according to $r$,
the lower bound of the part {\bf 3)} of Theorem \ref{t4} implies
the one of the part {\bf 2)} of the same Theorem. More precisely,
it can be checked that the necessary and sufficiently condition
for the holding of this improvement is $r \geq \frac{n^{\frac{1}{n
- 1}} + 1}{n^{\frac{1}{n - 1}} - 1}$, that is $n \leq f(r)$, where
$f$ is a real function which is equivalent to $\frac{1}{2} r
\log{r}$ as $r$ tends to infinite.
 \item[iv)] Under the additional assumptions $7 \leq r \leq 2
 u_0$ and $n \geq u_0 - \frac{3 r + 1}{2}$ (resp. $r \leq 2 u_0$ and $n \geq u_0 - \frac{3 r +
 1}{2}$), the lower bound of the part {\bf 4)} of Theorem
 \ref{t4} implies the one of the part {\bf 1)} (resp. {\bf 2)}) of
 the same Theorem up to the multiplicative constant $\frac{2}{\pi}$ (resp.
 $\frac{1}{\pi}$).\\
(Notice that the function $x \mapsto \sqrt{x} (x +
1)^{\frac{u_0}{x}}$ is decreasing on the interval $[7 , 2 u_0]$,
then if $7 \leq r \leq 2 u_0$, we have $\sqrt{r} (r +
1)^{\frac{u_0}{r}} > 2 u_0$).
 \item[v)] Now, we check that if $r \leq \frac{2}{3}
 u_0$ and $n \geq u_0 - \frac{3 r + 1}{2}$, the lower bound of the
 part {\bf 4)} of Theorem \ref{t4} implies (up to a multiplicative
 constant) the one of Conjecture \ref{conj1}. Indeed, if $r \leq \frac{2}{3}
 u_0$ and $n \geq u_0 - \frac{3 r + 1}{2}$, the decrease of the
 function $x \mapsto \sqrt{x} (x + 1)^{\frac{u_0}{x} - 1}$ on the
 interval $[1 , + \infty[$ implies:
 $\sqrt{r} (r + 1)^{\frac{u_0}{r} - 1} \geq \sqrt{\frac{2}{3} u_0} \left(\frac{2}{3} u_0 + 1\right)^{\frac{1}{2}}
 > \frac{2}{3} u_0$ which gives (by using the lower bound of the
 part {\bf 4)} of Theorem \ref{t4}):
$$\ppcm\{u_0 , \dots , u_n\} ~\geq~ \frac{2}{3 \pi} u_0 (r + 1)^n .$$
$\bullet$ More generally, for any given real $\xi \geq
 \frac{3}{2}$, if we suppose $r \leq \frac{1}{\xi} u_0$ and
 $n \geq u_0 - \frac{3 r + 1}{2}$ then the decrease of the
 function $x \mapsto \sqrt{x} (x + 1)^{\frac{u_0}{x} - \xi +
 \frac{1}{2}}$ on the interval $[1 , + \infty[$ implies:
$\sqrt{r} (r + 1)^{\frac{u_0}{r} - \xi + \frac{1}{2}} \geq
\sqrt{\frac{u_0}{\xi}} \left(\frac{u_0}{\xi} +
1\right)^{\frac{1}{2}} > \frac{u_0}{\xi}$ which gives (by using
the lower bound of the part {\bf 4)} of Theorem \ref{t4}):
$$\ppcm\{u_0 , \dots , u_n\} ~\geq~ \frac{1}{\pi \xi} u_0 (r + 1)^{n + \xi - \frac{3}{2}} .$$
Remark that if $\xi > \frac{3}{2}$, this lower bound is stronger
than the one of Conjecture \ref{conj1}.
\end{description}
\subsection{Results about the quadratic sequences}
\begin{thm}\label{t7}
Let $\mathbf{u} = {(u_k)}_{k \in \mathbb{N}}$ be a sequence of
integers whose general term has the form:
$$u_k ~=~ a k (k + t) + b ~~~~ (\forall~\!\! k \in \mathbb{N}) ,$$
with $(a , t , b) \in {\mathbb{Z}}^3$, $a \geq 1$, $t \geq 0$ and
$\pgcd\{a , b\} = 1$.\\
Also let $m$ and $n$ (with $m < n$) be two non-negative integers
for which none of the terms $u_k$ $(m \leq k \leq n)$ of
$\mathbf{u}$ is zero. Then the integer $\ppcm\{u_m , \dots ,
u_n\}$ is a multiple of the rational number:
$$A_{\mathbf{u}}(t , m , n) ~:=~ \begin{cases}2 \frac{u_0 \dots u_n}{(2 n)!} &
\text{if $(t , m) = (0 , 0)$}~\vspace{0.1cm} \\
(2 m + t - 1)! \frac{u_m \dots u_n}{(2 n + t)!} &
\text{otherwise}\end{cases} .$$
\end{thm}
\begin{coll}\label{c1}
Let $\mathbf{u} = {(u_k)}_{k \in \mathbb{N}}$ be a sequence of
integers as in the above Theorem and $n$ be a positive integer.
Then, if the $(n + 1)$ first terms $u_0 , \dots , u_n$ of the
sequence $\mathbf{u}$ are all non-zero, then we have:
$$\ppcm\{u_0 , \dots , u_n\} ~\geq~ \begin{cases}2 b \left(\frac{a}{4}\right)^n & \text{if $t = 0$}~\vspace{0.1cm} \\
\frac{b}{t 2^t} \left(\frac{a}{4}\right)^n & \text{if $t \geq
1$}\end{cases} .$$
\end{coll}
{\bf Remark.} It is clear that the lower bound of Corollary
\ref{c1} is nontrivial only if $a \geq 5$. Such as it is, this
corollary cannot thus give a nontrivial lower bound for the
numbers $\ppcm\{1^2 + 1 , 2^2 + 1 , \dots , n^2 + 1\}$ ($n \geq
1$). But we remark that if $r \geq 3$ is an integer, it gives a
nontrivial lower bound for the last common multiple of consecutive
terms of the sequence ${(r^2 n^2 + 1)}_{n \geq 1}$ which is a
subsequence of ${(n^2 + 1)}_n$. So we can first obviously bound
from below $\ppcm\{1^2 + 1 , 2^2 + 1 , \dots , n^2 + 1\}$ by
$\ppcm\{r^2 + 1 , r^2 2^2 + 1 , \dots , r^2 k^2 + 1\}$ (with $k :=
[\frac{n}{r}]$), then use Corollary \ref{c1} to bound from below
this new quantity. We obtain in this way:
$$
\ppcm\{1^2 + 1 , 2^2 + 1 , \dots , n^2 + 1\} ~\geq~ 2
\left(\frac{r^2}{4}\right)^k ~>~ 2
\left(\frac{r^2}{4}\right)^{\frac{n}{r} - 1} ~=~ \frac{8}{r^2}
\left\{\left(\frac{r}{2}\right)^{\frac{2}{r}}\right\}^n .
$$
This gives (for any choice of $r \geq 3$) a nontrivial lower bound
for the numbers $\ppcm\{1^2 + 1 , 2^2 + 1 , \dots , n^2 + 1\}$ $(n
\geq 1)$. We easily verify that the optimal lower bound
corresponds to $r = 5$, that is:
$$\ppcm\{1^2 + 1 , 2^2 + 1 , \dots , n^2 + 1\} ~\geq~ 0,32
(1,442)^n ~~~~ (\forall n \geq 1) .$$
\subsection{Results about the least common multiple of a finite number of consecutive integers}
\noindent The following Theorem is an immediate consequence of
Theorems \ref{t2} and \ref{t3}.
\begin{thm}\label{t8}
For any non-negative integer $k$ and any positive integer $n$,
the integer $\ppcm\{n , n + 1 , \dots , n + k\}$ is a multiple of the integer $n \binom{n + k}{k}$.\\
Further, if the congruence $n (n + k) \equiv 0 ~\!\mod (k!)$ is
satisfied, then we have precisely:
$$\ppcm\{n , n + 1 , \dots , n + k\} ~=~ n \binom{n + k}{k} .$$
\end{thm}
The following result is independent of all the results previously
quoted. It gives a multiple $m_{n , k}$ of the integer $\ppcm\{n ,
n + 1 , \dots , n + k\}$ $(k \in \mathbb{N} , n \in
{\mathbb{N}}^*)$ which is optimal and simple in its dependance on
$n$.
\begin{thm}\label{t9}
For any non-negative integer $k$ and any positive integer $n$, the
integer $\ppcm\{n , n + 1 , \dots , n + k\}$ divides the integer
$\displaystyle n \binom{n + k}{k} \ppcm\left\{\binom{k}{0}
, \binom{k}{1} , \dots , \binom{k}{k}\right\}$.~\vspace{0.2cm}\\
Further, if the congruence $n + k + 1 \equiv 0 ~\!\mod (k!)$ is
satisfied, then we have precisely:
$$\ppcm\{n , n + 1 , \dots , n + k\} ~=~ n \binom{n + k}{k} \ppcm\left\{\binom{k}{0} ,
\binom{k}{1} , \dots , \binom{k}{k}\right\} .$$
\end{thm}
\section{Proofs}
{\bf Proof of Lemma \ref{l1}.} In the situation of Lemma \ref{l1},
we have:
\begin{eqnarray*}
\frac{\ppcm\{\alpha_i , i \in I\} . \ppcm\{\beta_i , i \in
I\}}{\gamma} & = & \ppcm\{\alpha_i , i \in I\} . \ppcm\{\beta_i ,
i \in I\} \sum_{j \in I} \frac{1}{\alpha_j \beta_j} \\
& = & \sum_{j \in I} \frac{\ppcm\{\alpha_i , i \in I\}}{\alpha_j}
. \frac{\ppcm\{\beta_i , i \in I\}}{\beta_j} .
\end{eqnarray*}
This last sum is clearly an integer because for any $j \in I$, the
two numbers $\frac{\ppcm\{\alpha_i , i \in I\}}{\alpha_j}$ and
$\frac{\ppcm\{\beta_i , i \in I\}}{\beta_j}$ are integers. Lemma
\ref{l1} follows. $~~~~\blacksquare$ \vskip 1mm \noindent {\bf
Proof of Theorem \ref{t1}.} Theorem \ref{t1} follows by applying
Lemma \ref{l1} to the identity:
$$\sum_{j = 0}^{n} \frac{1}{u_j} . \frac{1}{\displaystyle \prod_{0 \leq i \leq n , i \neq j} \!\!\!(u_i - u_j)}
~=~ \frac{1}{u_0 u_1 \dots u_n} ,$$ which we obtain by taking $x =
0$ in the decomposition to simple elements of the rational
fraction $x \mapsto \frac{1}{(x + u_0) (x + u_1) \dots (x +
u_n)}$. $~~~~\blacksquare$ \vskip 1mm \noindent {\bf Proof of
Theorem \ref{t2}.} By replacing if necessary the sequence
${(u_n)}_n$ by the sequence with general term $v_n :=
\frac{u_n}{\pgcd\{u_0 , u_1\}}$ $(\forall n \in \mathbb{N})$, we
may assume that $u_0$ and $u_1$ are coprime. Under this
hypothesis, we have to show that the integer $\ppcm\{u_0 , \dots ,
u_n\}$ is a multiple of the rational number $\frac{u_0 \dots
u_n}{n!}$ (for any $n \in \mathbb{N}$).

Let $n$ be a fixed non-negative integer. From Theorem \ref{t1},
the integer $\ppcm\{u_0 , \dots , u_n\}$ is a multiple of the
rational number
$$\frac{u_0 \dots u_n}{\ppcm\left\{\displaystyle \prod_{0 \leq i \leq n , i \neq j} \!\!\!(u_i - u_j) ~;~
0 \leq j \leq n\right\}} .$$ Let $r$ denotes the difference of the
arithmetic sequence ${(u_k)}_k$. We have for any $(i , j) \in
{\mathbb{N}}^2$: $u_i - u_j = (i - j) r$, then for any $j \in \{0
, \dots , n\}$:
\begin{eqnarray*}
\prod_{0 \leq i \leq n , i \neq j} \!\!\!(u_i - u_j) & = &
\prod_{0 \leq i \leq n , i \neq j} \!\!\!(i - j) r \\
& = & r^n \left\{(-j) (1 - j) (2 - j) \dots (- 1)\right\} .
\left\{1 . 2 \dots (n - j)\right\} \\
& = & r^n (-1)^j j! (n - j)! .
\end{eqnarray*}
Hence:
\begin{eqnarray*}
\ppcm\left\{\displaystyle \prod_{0 \leq i \leq n , i \neq j}
\!\!\!(u_i - u_j) ~;~ 0 \leq j \leq n\right\} & = &
\ppcm\left\{r^n (-1)^j j! (n - j)! ~;~ 0 \leq j \leq n\right\}
\\
& = & r^n \ppcm\left\{j! (n - j)! ~;~ 0 \leq j \leq n\right\} \\
& = & r^n n!
\end{eqnarray*}
(because each integer $j! (n - j)!$ divides $n!$ and
for $j = 0$ or $n$, we have $j! (n - j)! = n!$).\\
Thus the integer $\ppcm\{u_0 , \dots , u_n\}$ is a multiple of the
rational number $\frac{u_0 \dots u_n}{r^n n!}$. But our hypothesis
``$u_0$ coprime with $u_1$''$~\!\!$ implies that $r$ is coprime
with all terms of the sequence ${(u_k)}_k$, which implies that
$r^n$ is coprime with the product $u_0 \dots u_n$. By the Gauss
lemma, we finally conclude that the integer $\ppcm\{u_0 , \dots ,
u_n\}$ is a multiple of the rational number $\frac{u_0 \dots
u_n}{n!}$ as required. $~~~~\blacksquare$ \vskip 1mm \noindent
{\bf Proof of Theorem \ref{t3}.} We need the following preliminary
Lemma: \vskip 1mm \noindent {\bf Lemma.---} {\it Let $n$ be a
positive integer and $x$ and $y$ be two integers satisfying: $x -
y \equiv 0 ~\mod(n)$ and $x y \equiv 0 ~\mod (n!)$. Then $x$ and
$y$ are multiples of
$n$.} \vskip 1mm \noindent {\bf Proof.} We distinguish the following four cases:\\
$\bullet$ \underline{If $n = 1$:} In this case, the result of Lemma is trivial.\\
$\bullet$ \underline{If $n$ is prime:} In this case, since $x^2 =
x (x - y) + x y$, we have $x^2 \equiv 0 ~\mod (n)$, but since $n$
is supposed prime, we conclude that $x \equiv 0 ~\mod (n)$ and then that $y = x - (x - y) \equiv 0 ~\mod (n)$.\\
$\bullet$ \underline{If $n = 4$:} In this case, we have $x - y
\equiv 0 ~\mod (4)$ and $x y \equiv 0 ~\mod (24)$ and we have to
show that $x$ and $y$ are multiples of $4$. Let us argue by
contradiction. Then, since $x \equiv y ~\mod (4)$, we have:\\
--- Either $x \equiv y \equiv 1 , 3 ~\mod (4)$ which implies $x y \equiv 1 ~\mod (4)$ and contradicts the congruence
$x y \equiv 0 ~\mod (24)$.\\
--- Or $x \equiv y \equiv 2 ~\mod (4)$ which implies $x y \equiv 4 ~\mod (8)$
and contradicts the congruence $x y \equiv 0 ~\mod (24)$ again.\\
Thus the Lemma holds for $n = 4$.\\
$\bullet$ \underline{If $n \geq 5$ and $n$ is not prime:} In this
case, it is easy to see that the integer $(n - 1)!$ is a multiple
of $n$, so that the integer $n!$ is a multiple of $n^2$. We thus
have $x - y \equiv 0 ~\mod (n)$ and $x y \equiv 0 ~\mod (n^2)$.

Let us argue by contradiction. Suppose that one at least of the
two integers $x$ and $y$ is not a multiple of $n$. To fix the
ideas, suppose for instance that $x \not\equiv 0 ~\mod (n)$. Then,
there exists a prime number $p$ dividing $n$ such that $v_p(x) <
v_p(n)$. But since $x y \equiv 0 ~\mod (n^2)$, we have $v_p(x y)
\geq v_p(n^2)$, that is $v_p(x) + v_p(y) \geq 2 v_p(n)$. This
implies that $v_p(y) \geq 2 v_p(n) - v_p(x) > v_p(x)$ (because
$v_p(x) < v_p(n)$). Thus, the $p$-adic valuations of the integers
$x$ and $y$ are distinct. Then we have: $v_p(x - y) = \min(v_p(x)
, v_p(y)) = v_p(x) < v_p(n)$, which contradicts the fact that $(x
- y)$ is a multiple of $n$. The Lemma is proved. ~\vspace{0.1cm}\\
{\bf Return to the proof of Theorem \ref{t3}:}~\vspace{0.1cm}

The case $n = 1$ is trivial. Next, we assume that $n \geq 2$. From
Theorem \ref{t2}, the integer $\ppcm\{u_0 , \dots , u_n\}$ is a
multiple of the rational number $\frac{u_0 \dots u_n}{n!}$. To
prove Theorem \ref{t3}, it remains to prove that $\frac{u_0 \dots
u_n}{n!}$ is also a multiple of $\ppcm\{u_0 , \dots , u_n\}$,
which means that $\frac{u_0 \dots u_n}{n!}$ is a multiple of each
of integers $u_0 , \dots , u_n$. Since $u_0 u_n$ is assumed a
multiple of $n!$, the number $\frac{u_0 \dots u_n}{n!}$ is
obviously a multiple of each of integers $u_1 , \dots , u_{n -
1}$. To conclude, it only remains to prove that this same number
$\frac{u_0 \dots u_n}{n!}$ is a multiple of $u_0$ and $u_n$, which
is equivalent to prove that the two integers $u_1 \dots u_n$ and
$u_0 \dots u_{n - 1}$ are multiples of $n!$. We first prove that
$u_0$ and $u_n$ are multiples of $n$. Denoting $r$ the difference
of the arithmetic sequence ${(u_k)}_k$, we have $u_n - u_0 = r n
\equiv 0 ~\mod (n)$ and $u_0 u_n \equiv 0 ~\mod (n!)$ (by
hypothesis). This implies (from the above Lemma) that $u_0$ and
$u_n$ effectively are multiples of $n$.

We now prove that the two integers $u_1 \dots u_n$ and $u_0 \dots
u_{n - 1}$ are multiples of $n!$. For any $1 \leq k \leq n - 1$,
we have: $u_k = u_0 + k r \equiv k r ~\mod (u_0)$, then:
$$
u_1 \dots u_{n - 1} ~\equiv~ (1 . r) (2 . r) \dots ((n - 1) . r)
~\mod (u_0) ~\equiv~ (n - 1)! r^{n - 1} ~\mod (u_0) .
$$
It follows that:
$$u_1 \dots u_{n - 1} u_n ~\equiv~ (n - 1)! u_n r^{n - 1} ~\mod (u_0 u_n) .$$
Since $u_n$ is a multiple of $n$ and (by hypothesis) $u_0 u_n$ is
a multiple of $n!$, the last congruence implies that $u_1 \dots
u_{n - 1} u_n$ is a multiple of $n!$.

Similarly, for any $1 \leq k \leq n - 1$, we have: $u_{n - k} =
u_n - k r \equiv - k r ~\mod (u_n)$, then:
$$
u_{n - 1} \dots u_1 ~\equiv~ (- (n - 1) . r) \dots (- 1 . r) ~\mod
(u_n) ~\equiv~ (- 1)^{n - 1} (n - 1)! r^{n - 1} ~\mod (u_n) .
$$
It follows that:
$$u_0 u_1 \dots u_{n - 1} ~\equiv~ (- 1)^{n - 1} (n - 1)! u_0 r^{n - 1} ~\mod (u_0 u_n) .$$
Since $u_0$ is a multiple of $n$ and (by hypothesis) $u_0 u_n$ is
a multiple of $n!$, the last congruence implies that $u_0 \dots
u_{n - 1}$ is also a multiple of $n!$. This completes the proof of
Theorem \ref{t3}. $~~~~\blacksquare$ \vskip 1mm \noindent {\bf
Proof of Theorem \ref{t4}.} For any integer $k \in \{0 , \dots ,
n\}$, the integer $\ppcm\{u_0 , \dots , u_n\}$ is obviously a
multiple of the integer $\ppcm\{u_k , \dots , u_n\}$ and from
Theorem \ref{t2}, this last integer is a multiple of the rational
number $\frac{u_k \dots u_n}{(n - k)!}$. It follows that for any
$k \in \{0 , \dots , n\}$, we have:
\begin{equation}\label{eq1}
\ppcm\{u_0 , \dots , u_n\} ~\geq~ \frac{u_k \dots u_n}{(n - k)!} .
\end{equation}
The idea consists in choosing $k$ as a function of $n , r$ and
$u_0$ in order to optimize the lower bound (\ref{eq1}), that is to
make the quantity $\frac{u_k \dots u_n}{(n - k)!}$ maximal.

Let ${(v_k)}_{0 \leq k \leq n}$ denotes the finite sequence of
general term: $v_k := \frac{u_k \dots u_n}{(n -
k)!}$. We have the following intermediate Lemma:~\vspace{0.1cm}\\
{\bf Lemma.---} {\it The sequence ${(v_k)}_{0 \leq k \leq n}$
reaches its maximum value at $$k_0 := \max\left\{0 ,
\left[\frac{n - u_0}{r + 1}\right] + 1\right\} .$$}~\vspace{0.1cm}\\
{\bf Proof.} For any $k \in \{0 , \dots , n - 1\}$, we have:
$\frac{v_{k + 1}}{v_k} = \frac{u_{k + 1} \dots u_n}{(n - k - 1)!}
/ \frac{u_k \dots u_n}{(n - k)!} = \frac{n - k}{u_k} = \frac{n -
k}{u_0 + k r}$, hence:
$$
v_{k + 1} \geq v_k ~\Longleftrightarrow~ \frac{n - k}{u_0 + k r}
\geq 1 ~\Longleftrightarrow~ k \leq \frac{n - u_0}{r + 1}
~\Longleftrightarrow~ k \leq \left[\frac{n - u_0}{r + 1}\right] .
$$
This permits us to determine the variations of the finite sequence
${(v_k)}_{0 \leq k \leq n}$ according to the position of $n$
compared to $u_0$. If $n < u_0$, the sequence ${(v_k)}_{0 \leq k
\leq n}$ is decreasing and it thus reaches its maximum value at $k
= 0$. In the other case i.e $n \geq u_0$, the sequence ${(v_k)}_{0
\leq k \leq n}$ is increasing until the integer $\left[\frac{n -
u_0}{r + 1}\right] + 1$ then it decreases, so it reaches its
maximum value at $k = \left[\frac{n - u_0}{r + 1}\right] + 1$. The
Lemma follows.~\vspace{0.1cm}

The following intermediary lemma gives an identity which permits
to bound from bellow $v_k$ by simple expressions (as function as
$u_0 , r$ and $n$) for the particular values of $k$ which are
rather close to the integer $k_0$ of the above Lemma.~\vspace{0.1cm}\\
{\bf Lemma.---} {\it For any $k \in \{0 , \dots , n\}$, we have:
\begin{equation}\label{eq2}
v_k ~=~ \frac{r^{n - k + 1}}{\int_{0}^{1} x^{k + \frac{u_0}{r} -
1} (1 - x)^{n - k} d x} .
\end{equation}
}~\vspace{0.1cm}\\
{\bf Proof.} For any $0 \leq k \leq n$, we have:
\begin{eqnarray*}
v_k ~:=~ \frac{u_k \dots u_n}{(n - k)!} & = & \frac{u_k \left(u_k
+ r\right) \dots \left(u_k + (n - k) r\right)}{(n - k)!} \\
& = & r^{n - k + 1} \frac{\frac{u_k}{r} \left(\frac{u_k}{r} +
1\right) \dots \left(\frac{u_k}{r} + n - k\right)}{(n - k)!} \\
& = & r^{n - k + 1} \frac{\Gamma\left(\frac{u_k}{r} + n - k +
1\right)}{\Gamma\left(\frac{u_k}{r}\right) . \Gamma\left(n - k +
1\right)} \\
& = & \frac{r^{n - k + 1}}{\beta\left(\frac{u_k}{r} , n - k +
1\right)} ,
\end{eqnarray*}
where $\Gamma$ and $\beta$ denote the Euler's functions. The
identity (\ref{eq2}) of Lemma follows from the well known integral
formula of the $\beta$-function. The Lemma is
proved.~\vspace{0.1cm}

Because of some technical difficulties concerning the lower bound
of the right-hand side of (\ref{eq2}) for $k = k_0$, we are led to
bound from below this side for other values of $k$ which are close
to $k_0$. So, we obtain the lower bounds of the parts {\bf 1)} and
{\bf 4)} of Theorem \ref{t4} by bounding from below $v_k$ for $k =
\left[\frac{n - 1}{r + 1} + 1\right]$ and for the nearest integer
$k$ to the real $\frac{n + r - u_0}{r + 1}$ respectively. Further,
we obtain the remaining parts {\bf 2)} and {\bf 3)} of Theorem
\ref{t4} by another method which doesn't use the identity
(\ref{eq2}). We first prove the parts {\bf 1)} and {\bf 4)} of Theorem \ref{t4}.~\vspace{0.1cm}\\
{\bf Proof of the part 1) of Theorem \ref{t4}:}~\vspace{0.1cm}

Let $k_1 := \left[\frac{n - 1}{r + 1} + 1\right]$. Using the
identity (\ref{eq2}), we are going to get a lower bound for
$v_{k_1}$ which depends on $u_0 , r$ and $n$. The integer $k_1$
satisfies $\frac{n - 1}{r + 1} < k_1 \leq \frac{n - 1}{r + 1} + 1
= \frac{n + r}{r + 1}$. We thus have:
\begin{equation}\label{eq3}
r^{n - k_1 + 1} ~\geq~ r^{\frac{(n - 1) r}{r + 1} + 1}
\end{equation}
and for any real $x \in [0 , 1]$:
$$x^{k_1 + \frac{u_0}{r} - 1} (1 - x)^{n - k_1} ~\leq~ x^{\frac{n - 1}{r + 1} + \frac{u_0}{r} - 1}
(1 - x)^{\frac{(n - 1) r}{r + 1}} ,$$ which gives:
\begin{equation}\label{eq4}
\int_{0}^{1} x^{k_1 + \frac{u_0}{r} - 1} (1 - x)^{n - k_1} d x
~\leq~ \int_{0}^{1} {\left\{x (1 - x)^r\right\}}^{\frac{n - 1}{r +
1}} x^{\frac{u_0}{r} - 1} d x .
\end{equation}
By studying the function $x \mapsto x (1 - x)^r$, we may show that
for any real $x \in [0 , 1]$, we have: $x (1 - x)^r \leq
\frac{r^r}{(r + 1)^{r + 1}}$. Substituting this into the
right-hand side of (\ref{eq4}), we deduce that:
\begin{equation}\label{eq5}
\int_{0}^{1} x^{k_1 + \frac{u_0}{r} - 1} (1 - x)^{n - k_1} d x
~\leq~ \frac{r^{\frac{(n - 1) r}{r + 1}}}{(r + 1)^{n - 1}} .
\frac{r}{u_0}
\end{equation}
By combining the two relations (\ref{eq3}) and (\ref{eq5}), we
finally obtain:
$$\frac{r^{n - k_1 + 1}}{\int_{0}^{1} x^{k_1 + \frac{u_0}{r} - 1} (1 - x)^{n - k_1} d x} ~\geq~ u_0 (r + 1)^{n - 1} .$$
Then the first lower bound of the part 1) of Theorem \ref{t4}
follows from the relations (\ref{eq2}) and (\ref{eq1}).

If $n$ is a multiple of $(r + 1)$, the second lower bound of the
part 1) of Theorem \ref{t4} follows by taking in the above proof
instead of $k_1$ the integer $k = \frac{n}{r + 1}$.~\vspace{1mm} \\
{\bf Proof of the part 4) of Theorem \ref{t4}:}~\vspace{1mm}

The particular case $n = 0$ of the part 4) of Theorem \ref{t4}
follows from the fact that the function $x \mapsto \sqrt{x} (x +
1)^{\frac{u_0}{x} - 1}$ is decreasing on the interval $[1 , +
\infty[$. Next, we suppose that $n \geq 1$. The hypothesis $n \geq
u_0 - \frac{3 r + 1}{2}$ means that the real $\frac{n + r - u_0}{r
+ 1}$ is greater than or equal to $- \frac{1}{2}$. Since this same
real $\frac{n + r - u_0}{r + 1}$ is less than or equal to $n +
\frac{1}{2}$ (because $n \geq 1$), then there exists an integer
$k_2 \in \{0 , \dots , n\}$ satisfying:
$$- \frac{1}{2} ~\leq~ k_2 - \frac{n + r - u_0}{r + 1} ~\leq~ \frac{1}{2} .$$
It follows that:
\begin{equation}\label{eq6}
r^{n - k_2 + 1} ~\geq~ r^{\frac{r (n - 1) + u_0}{r + 1} +
\frac{1}{2}}
\end{equation}
and that for any real $x \in ]0 , 1[$:
\begin{eqnarray*}
x^{k_2 + \frac{u_0}{r} - 1} (1 - x)^{n - k_2} & \leq & x^{\frac{r
(n - 1) + u_0}{r (r + 1)} - \frac{1}{2}} (1 - x)^{\frac{r (n - 1)
+ u_0}{r + 1} - \frac{1}{2}} \\
& = & \left\{x (1 - x)^r\right\}^{\frac{r (n - 1) + u_0}{r (r +
1)}} \frac{1}{\sqrt{x (1 - x)}} \\
& \leq & \left(\frac{r^r}{(r + 1)^{(r + 1)}}\right)^{\frac{r (n -
1) + u_0}{r (r + 1)}} \frac{1}{\sqrt{x (1 - x)}}
\end{eqnarray*}
(because $x (1 - x)^r \leq \frac{r^r}{(r + 1)^{r + 1}}$ for any $x
\in [0 , 1]$).\\
Consequently:
$$\int_{0}^{1} x^{k_2 + \frac{u_0}{r} - 1} (1 - x)^{n - k_2} d x ~\leq~
\left(\frac{r^r}{(r + 1)^{(r + 1)}}\right)^{\frac{r (n - 1) +
u_0}{r (r + 1)}} \int_{0}^{1} \frac{d x}{\sqrt{x (1 - x)}} .$$
Since $\int_{0}^{1} \frac{d x}{\sqrt{x (1 - x)}} = \pi$, we deduce
that:
\begin{equation}\label{eq7}
\int_{0}^{1} x^{k_2 + \frac{u_0}{r} - 1} (1 - x)^{n - k_2} d x
~\leq~ \pi \frac{r^{\frac{r (n - 1) + u_0}{r + 1}}}{(r + 1)^{n - 1
+ \frac{u_0}{r}}} .
\end{equation}
By combining the two relations (\ref{eq6}) and (\ref{eq7}), we
finally obtain:
$$\frac{r^{n - k_2 + 1}}{\int_{0}^{1} x^{k_2 + \frac{u_0}{r} - 1} (1 - x)^{n - k_2}} ~\geq~ \frac{1}{\pi}
\sqrt{r} (r + 1)^{n - 1 + \frac{u_0}{r}}$$ and we conclude the
lower bound of the part 4) of Theorem \ref{t4} by using the
identity (\ref{eq2}) and the lower bound (\ref{eq1}).~\vspace{0.1cm}\\

We obtain the two remaining parts 2) and 3) of Theorem \ref{t4} by
using the same idea which consists to bound from below $v_k =
\frac{u_k \dots u_n}{(n - k)!}$ for some particular values of $k
\in \{0 , \dots , n\}$. The only difference with the last parts 1)
and 4) proved above is that here such particular values are not
explicit,
we just show their existence by using the following Lemma:~\vspace{0.1cm}\\
{\bf Lemma.---} Let $x$ be a real and $n$ be a positive integer.
Then:
\begin{description}
\item[1)] there exists an integer $k$ $(1 \leq k \leq n)$ such
that:
$$k \binom{n}{k} x^{n - k + 1} ~\geq~ x (x + 1)^{n - 1} .$$
\item[2)] There exists an odd integer $\ell$ $(1 \leq \ell \leq
n)$ such that:
$$\ell \binom{n}{\ell} x^{n - \ell + 1} ~\geq~ \frac{n}{n + 1} x \left\{(x + 1)^{n - 1} +
(x - 1)^{n - 1}\right\} .$$
\end{description}~\vspace{0.1cm}\\
{\bf Proof.} The first part of Lemma follows from the identity:
\begin{equation}\label{eq8}
\sum_{k = 1}^{n} k \binom{n}{k} x^{n - k + 1} ~=~ n x (x + 1)^{n -
1}
\end{equation}
which can be proved by deriving with respect to $u$ the binomial
formula\\ $\sum_{k = 0}^{n} \binom{n}{k} u^k x^{n - k} = (u +
x)^n$ and then by taking $u = 1$ in the obtained formula.\\
The second part of Lemma follows from the identity:
\begin{equation}\label{eq9}
\sum_{\begin{array}{c}\scriptstyle{1 \leq k \leq n} \\
\scriptstyle{k ~\!\text{odd}}\end{array}} \!\!\!k \binom{n}{k}
x^{n - k + 1} ~=~ \frac{1}{2} n x \left\{(x + 1)^{n - 1} + (x -
1)^{n - 1}\right\}
\end{equation}
which follows from (\ref{eq8}) by remarking that:
$$\sum_{\begin{array}{c}\scriptstyle{1 \leq k \leq n} \\
\scriptstyle{k ~\!\text{odd}}\end{array}} \!\!\!k \binom{n}{k}
x^{n - k + 1} ~=~ \frac{1}{2} \left\{\sum_{k = 1}^{n} k
\binom{n}{k} x^{n - k + 1} + (- 1)^n \sum_{k = 1}^{n} k
\binom{n}{k} (- x)^{n - k + 1}\right\} .$$ The Lemma is proved.~\vspace{0.1cm}\\
{\bf Proof of the parts 2) and 3) of Theorem
\ref{t4}:}~\vspace{0.1cm}

We have for any $k \in \{1 , \dots , n\}$:
$$v_k := \frac{u_k \dots u_n}{(n - k)!} ~\geq~ \frac{(k r) ((k + 1) r) \dots (n r)}{(n - k)!} = k \binom{n}{k}
r^{n - k + 1} .$$ These lower bounds of $v_k$ $(1 \leq k \leq n)$
implie (by using the above Lemma) that there exist an integer $k
\in \{1 , \dots , n\}$ and an odd integer $\ell \in \{1 , \dots ,
n\}$ for which we have:
$$v_k ~\geq~ r (r + 1)^{n - 1} ~~\text{and}~~ v_{\ell} ~\geq~ \frac{n}{n + 1} r \left\{(r + 1)^{n - 1}
+ (r - 1)^{n - 1}\right\} .$$ We conclude by using the relation
(\ref{eq1}). This completes the proof of Theorem \ref{t4}.
$~~\blacksquare$ \vskip 1mm \noindent {\bf Proof of Theorem
\ref{t5}.} Let us argue by contradiction. Then, we can find a
rational number $\frac{a}{b} > \frac{3}{2}$ (with $a , b$ are
positive integers) for which we have for any arithmetic
progression ${(u_k)}_k$ with positive difference $r$, satisfying
the hypothesis of Theorem \ref{t4} and for any non-negative
integer $n \geq u_0 - \frac{a}{b} r - \frac{1}{2}$:
$$\ppcm\{u_0 , \dots , u_n\} ~\geq~ \frac{1}{\pi} \sqrt{r} (r + 1)^{n - 1 + \frac{u_0}{r}} .$$
We introduce a non-negative parameter $\delta$ and the arithmetic
progression ${(u_k)}_k$ (depending on $\delta$) with first term
$u_0 := a b \delta + 1$ and difference $r := b^2 \delta$. The
integers $u_0$ and $r$ are coprime because they verify the Bézout
identity $(1 - a b \delta) u_0 + a^2 \delta r = 1$. The sequence
${(u_k)}_k$ thus satisfies all the hypotheses of Theorem \ref{t4}.
Since the integer $n = 1$ satisfies $n \geq u_0 - \frac{a}{b} r -
\frac{1}{2} = \frac{1}{2}$, we must have:
\begin{equation}\label{eq10}
\ppcm\{u_0 , u_1\} ~\geq~ \frac{1}{\pi} \sqrt{r} (r +
1)^{\frac{u_0}{r}} .
\end{equation}
Further, we have
$$\ppcm\{u_0 , u_1\} = u_0 u_1 = (a b \delta + 1) \left((a b + b^2) \delta + 1\right) = O(\delta^2)$$
and
$$\frac{1}{\pi} \sqrt{r} (r + 1)^{\frac{u_0}{r}} = \frac{1}{\pi} b \sqrt{\delta} (b^2 \delta + 1)^{\frac{a}{b} +
\frac{1}{b^2 \delta}} = O\left(\delta^{\frac{a}{b} +
\frac{1}{2}}\right) .$$ But since $\frac{a}{b} + \frac{1}{2} > 2$,
The relation (\ref{eq10}) cannot holds for $\delta$ sufficiently
large. Contradiction. Theorem \ref{t5} follows. $~~~~\blacksquare$
\vskip 1mm \noindent {\bf Proof of Theorem \ref{t6}.} Let us prove
the assertion 1) of Theorem \ref{t6}. The fact that the constant
$C$ of this assertion is greater than or equal to $\frac{1}{\pi}$
is an immediate consequence of the part 4) of Theorem \ref{t4}. In
order to prove the upper bound $C \leq \frac{3}{2}$, we introduce
a parameter $\delta \in \mathbb{N}$ and the arithmetic sequence
${(u_k)}_k$ (depending on $\delta$), with first term $u_0 := 3
\delta + 2$ and difference $r := 2 \delta + 1$. The integers $u_0$
and $r$ are coprime because they verify the Bézout identity $2 u_0
- 3 r = 1$. So, this sequence ${(u_k)}_k$ satisfies all the
hypotheses of Theorem \ref{t4}. Since $u_0 - \frac{3 r + 1}{2} =
0$, we must have for any non-negative integer $n$: $\ppcm\{u_0 ,
\dots , u_n\} \geq C \sqrt{r} (r + 1)^{n - 1 + \frac{u_0}{r}}$, in
particular (for $n = 0$): $u_0 \geq C \sqrt{r} (r +
1)^{\frac{u_0}{r} - 1}$, hence:
$$C ~\leq~ \frac{u_0}{\sqrt{r} (r + 1)^{\frac{u_0}{r} - 1}} .$$
Since this last upper bound holds for any $\delta \in \mathbb{N}$,
we finally deduce that:
$$C ~\leq~ \lim_{\delta \rightarrow + \infty} \frac{u_0}{\sqrt{r} (r + 1)^{\frac{u_0}{r} - 1}}
~=~ \lim_{\delta \rightarrow + \infty}
\frac{3 \delta + 2}{\sqrt{2 \delta + 1} (2 \delta +
2)^{\frac{\delta + 1}{2 \delta + 1}}} ~=~ \frac{3}{2}$$ as
required.

Now, let us prove the assertion 2) of Theorem \ref{t6}. Let $n_0$
be a fixed non-negative integer. As above, the lower bound $C(n_0)
\geq \frac{1}{\pi}$ is an immediate consequence of the part 4) of
Theorem \ref{t4}. In order to prove the upper bound of Theorem
\ref{t6} for the constant $C(n_0)$, we choose an integer $n_1$
such that $n_0 + 3 \leq n_1 \leq 2 n_0 + 6$ and that $(n_1 + 1)$
is prime (this is possible from the Bertrand postulate). Then, we
introduce a parameter $\delta \in \mathbb{N}$ which is not a
multiple of $(n_1 + 1)$ and the arithmetic progression ${(u_k)}_k$
(depending on $\delta$), with first term $u_0 := 3 \delta n_1!$
and difference $r := 2 \delta n_1! + n_1 + 1$. These integers
$u_0$ and $r$ are coprime. Indeed, a common divisor $d \geq 1$
between $u_0$ and $r$ divides $3 r - 2 u_0 = 3 (n_1 + 1)$, thus it
divides $\pgcd\{u_0 , 3 (n_1 + 1)\} = 3 \pgcd\{\delta n_1! , n_1 +
1\}$. Further, the fact that $(n_1 + 1)$ is prime implies that
$(n_1 + 1)$ is coprime with $n_1!$, moreover since $\delta$ is not
a multiple of $(n_1 + 1)$, the integer $(n_1 + 1)$ also is coprime
with $\delta$. It follows that $(n_1 + 1)$ is coprime with the
product $\delta n_1!$. Hence $d$ divides $3$. But since $3$
divides $2 \delta n_1!$ (because $n_1 \geq 3$) and $3$ doesn't
divide $n_1 + 1$ (because $n_1 + 1$ is a prime number $\geq 5$)
then $3$ cannot divide the sum $2 \delta n_1! + (n_1 + 1) = r$,
which proves that $d \neq 3$. Consequently $d = 1$, that is $u_0$
and $r$ are coprime effectively. The sequence ${(u_k)}_k$ which we
have introduced thus satisfies all the hypotheses of Theorem
\ref{t4}. Since $n_1 \geq \max\{n_0 , u_0 - \frac{3 r + 1}{2}\}$
(because $n_1 \geq n_0 + 3$ and $u_0 - \frac{3 r + 1}{2} = -
\frac{3}{2} n_1 - 2 < 0$), then we must have $\ppcm\{u_0 , \dots ,
u_{n_1}\} \geq C(n_0) \sqrt{r} (r + 1)^{n_1 - 1 + \frac{u_0}{r}}$.
This gives:
$$C(n_0) ~\leq~ \frac{\ppcm\{u_0 , \dots , u_{n_1}\}}{\sqrt{r} (r + 1)^{n_1 - 1 + \frac{u_0}{r}}} .$$
Now, since $u_0$ is a multiple of $n_1!$, we have from Theorem
\ref{t3}: $\ppcm\{u_0 , \dots , u_{n_1}\} = \frac{u_0 \dots
u_{n_1}}{n_1!}$. Hence:
$$C(n_0) ~\leq~ \frac{u_0 \dots u_{n_1}}{n_1! \sqrt{r} (r + 1)^{n_1 - 1 + \frac{u_0}{r}}} .$$
Since this last upper bound of $C(n_0)$ holds for any $\delta \in
\mathbb{N}$ which is not a multiple of $(n_1 + 1)$, then we deduce
that:
\begin{equation}\label{eq11}
C(n_0) ~\leq \lim_{\begin{array}{c}\scriptstyle{\delta \rightarrow
+ \infty}
\\ \scriptstyle{\delta \not\equiv 0 ~\mod (n_1 + 1)}\end{array}}\!\!\!\!\!\frac{u_0 \dots u_{n_1}}{n_1!
\sqrt{r} (r + 1)^{n_1 - 1 + \frac{u_0}{r}}} .
\end{equation}
Let us calculate the limit from the right-hand side of
(\ref{eq11}). We have:
\begin{equation*}
\begin{split}
u_0 \dots u_{n_1} = \prod_{k = 0}^{n_1} (u_0 + k r) &= \prod_{k =
0}^{n_1} \left\{(2 k + 3) n_1! \delta + k (n_1 + 1)\right\} \\
&\sim_{+ \infty} \left({n_1!}^{n_1 + 1} \prod_{k = 0}^{n_1} (2 k +
3)\right) \delta^{n_1 + 1}
\end{split}
\end{equation*}
and:
\begin{equation*}
\begin{split}
\sqrt{r} (r + 1)^{n_1 - 1 + \frac{u_0}{r}} &= \left(2 \delta n_1!
+ n_1 + 1\right)^{1/2} \left(2 \delta n_1! + n_1 + 2\right)^{n_1 -
1 + \frac{3 \delta n_1!}{2 \delta n_1! + n_1 + 1}} \\
&\sim_{+ \infty} (2 \delta n_1!)^{n_1 + 1} .
\end{split}
\end{equation*}
Then:
$$
\frac{u_0 \dots u_{n_1}}{n_1! \sqrt{r} (r + 1)^{n_1 - 1 +
\frac{u_0}{r}}} ~\sim_{+ \infty}~ \frac{\displaystyle \prod_{k =
0}^{n_1} (2 k + 3)}{2^{n_1 + 1} n_1!} = \frac{(n_1 + 1) (n_1 +
\frac{3}{2})}{4^{n_1}} \binom{2 n_1 + 1}{n_1} .
$$
In the other words:
$$\lim_{\delta \rightarrow + \infty} \frac{u_0 \dots u_{n_1}}{n_1! \sqrt{r} (r + 1)^{n_1 - 1 +
\frac{u_0}{r}}} ~=~ \frac{(n_1 + 1) (n_1 + \frac{3}{2})}{4^{n_1}}
\binom{2 n_1 + 1}{n_1} .$$ It is easy to show (by induction on
$k$) that for any non-negative integer $k$, we have $\binom{2 k +
1}{k} < \sqrt{2} \frac{4^k}{\sqrt{k + \frac{3}{2}}}$. Using this
estimate for $k = n_1$, we finally deduce that:
\begin{eqnarray*}
\lim_{\delta \rightarrow + \infty} \frac{u_0 \dots u_{n_1}}{n_1!
\sqrt{r} (r + 1)^{n_1 - 1 + \frac{u_0}{r}}} & < & \sqrt{2} (n_1 +
1) \sqrt{n_1 + \frac{3}{2}} \\
& < &  4 (n_0 + 4) \sqrt{n_0 + 4} ~~~~ \text{(because $n_1 \leq 2
n_0 + 6$)} .
\end{eqnarray*}
The upper bound $C(n_0) < 4 (n_0 + 4) \sqrt{n_0 + 4}$ follows by
substituting this last estimate into (\ref{eq11}). This completes
the proof of Theorem \ref{t6}. $~~~~\blacksquare$ \vskip 1mm
\noindent {\bf Proof of Theorem \ref{t7}.} We first prove Theorem
\ref{t7} in the particular case $m = 0$. We deduce the general
case of the same Theorem by shifting the terms of the sequence
$\mathbf{u} = {(u_k)}_k$. Let $\mathbf{u}$ be a
sequence as in Theorem \ref{t7}.~\vspace{0.1cm}\\
$\bullet$ \underline{The case $m = 0$:} From Theorem \ref{t1}, the
integer $\ppcm\{u_0 , \dots , u_n\}$ is a multiple of the rational
number
\begin{equation}\label{eq*}
R ~:=~ \frac{u_0 \dots u_n}{\ppcm\left\{\displaystyle \prod_{0
\leq i \leq n , i \neq j} \!\!\!(u_i - u_j) ~;~ j = 0 , \dots ,
n\right\}} .
\end{equation}
Now, since we have for any $i , j \in \mathbb{N}$:
$$
u_i - u_j = \{a i (i + t) + b\} - \{a j (j + t) + b\} = a (i - j)
(i + j + t) ,
$$
then:
\begin{eqnarray*}
\prod_{0 \leq i \leq n , i \neq j} \!\!\!(u_i - u_j) & = &
\prod_{0 \leq i \leq n , i \neq j} \!\!\!\{a (i - j) (i + j + t)\} \\
& = & a^n \prod_{0 \leq i \leq n , i \neq j} \!\!\!(i - j) .
\prod_{0 \leq i \leq n , i \neq j} \!\!\!(i + j + t) \\
& = & \begin{cases} a^n (- 1)^j \frac{(n - j)! (n + j)!}{2} &
\text{if $t = 0$} ~\vspace{0.1cm}\\ a^n (- 1)^j \frac{(n -j)! (n +
j + t)!}{\varphi(j , t)} \frac{1}{2 j + t} & \text{if $t \geq
1$}\end{cases} ,
\end{eqnarray*}
where $\varphi(j , t) := 1$ if $t = 1$ and $\varphi(j , t) := (j +
1) \dots (j + t - 1)$ if $t \geq 2$. Since $(n - j)! (n + j + t)!$
divides $(2 n + t)!$ (because $\frac{(2 n + t)!}{(n - j)! (n + j +
t)!} = \binom{2 n + t}{n - j} \in \mathbb{N}$) and (if $t \geq 1$)
the integer $\varphi(j , t)$ is a multiple of $(t - 1)!$ (because
$\frac{\varphi(j , t)}{(t - 1)!} = \binom{j + t - 1}{t - 1} \in
\mathbb{N}$), then the product $\displaystyle \prod_{0 \leq i \leq
n , i \neq j} \!\!\!\!\!(u_i - u_j)$ divides the integer (which
does not depend on $j$):
$$f(t , n) ~:=~ \begin{cases}a^n \frac{(2 n)!}{2} & \text{if $t = 0$} ~\vspace{0.1cm}\\
a^n \frac{(2 n + t)!}{(t - 1)!} & \text{if $t \geq 1$}\end{cases}
.$$ Since $j$ is arbitrary in $\{0 , \dots , n\}$, then the
integer $\ppcm\{ \prod_{0 \leq i \leq n , i \neq j} (u_i - u_j) ;
j = 0 , \dots , n\}$ divides the integer $f(t , n)$. It follows
that the rational number $R$ (of (\ref{eq*})) is a multiple of the
rational number $\frac{u_0 \dots u_n}{f(t , n)} =
\frac{A_{\mathbf{u}}(t , 0 , n)}{a^n}$. Consequently, the integer
$\ppcm\{u_0 , \dots , u_n\}$ is a multiple of the rational number
$\frac{A_{\mathbf{u}}(t , 0 , n)}{a^n}$. Finally, since each term
of the sequence $\mathbf{u}$ is coprime with $a$ (because
$\pgcd\{a , b\} = 1$), we conclude from the Gauss Lemma that the
integer $\ppcm\{u_0 , \dots , u_n\}$ is a multiple of the rational
number $A_{\mathbf{u}}(t , 0 , n)$ as required.~\vspace{0.1cm}\\
$\bullet$ \underline{The general case ($m \in \mathbb{N}$):} Let
us consider the new sequence $\mathbf{v} = {(v_k)}_{k \in
{\mathbb{N}}}$ with general term:
$$
v_k := u_{k + m} = a' k (k + t') + b' ,
$$
where $a' := a$, $t' := 2 m + t$ and $b' := a m (m + t) + b$.\\
Since these integers $a'$, $t'$ and $b'$ verify $a' \geq 1$, $t'
\geq 0$ and $\pgcd\{a' , b'\} = \pgcd\{a , b\} = 1$ obviously,
then the sequence $\mathbf{v}$ satisfies all the hypotheses of
Theorem \ref{t7}. Thus, from the particular case (proved above) of
this Theorem, the integer $\ppcm\{v_0 , \dots , v_{n - m}\} =
\ppcm\{u_m , \dots , u_n\}$ is a multiple of the rational number
$A_{\mathbf{v}}(t' , 0 , n - m) = A_{\mathbf{u}}(t , m , n)$ which
provides the desired conclusion. $~~~~\blacksquare$ \vskip 1mm
\noindent {\bf Proof of Corollary \ref{c1}.} From Theorem
\ref{t7}, the integer $\ppcm\{u_0 , \dots , u_n\}$ is a multiple
of the rational number:
$$A_{\mathbf{u}}(t , 0 , n) ~:=~ \begin{cases} 2 \frac{u_0 \dots u_n}{(2 n)!} & \text{if $t = 0$} ~\vspace{0.1cm}\\
(t - 1)! \frac{u_0 \dots u_n}{(2 n + t)!} & \text{if $t \geq
1$}\end{cases} .$$ Let us get a lower bound for this last number
which doesn't depend on the terms of the sequence $\mathbf{u}$.
Using the obvious lower bounds $u_k \geq a k (k + t)$ $(1 \leq k
\leq n)$, we have:
$$
u_0 \dots u_n ~\geq~ b \{a . 1 . (1 + t)\} \{a . 2 . (2 + t)\}
\dots \{a . n . (n + t)\} = b a^n \frac{n! (n + t)!}{t!} ,
$$
then:
$$
A_{\mathbf{u}}(t , 0 , n) ~\geq~ \begin{cases}2 b
\frac{a^n}{\binom{2 n}{n}} & \text{if $t = 0$} ~\vspace{0.1cm}\\
\frac{b}{t} \frac{a^n}{\binom{2 n + t}{n}} & \text{if $t \geq
1$}\end{cases} ~\geq~ \begin{cases}2 b
\left(\frac{a}{4}\right)^n & \text{if $t = 0$} ~\vspace{0.1cm}\\
\frac{b}{t 2^t} \left(\frac{a}{4}\right)^n & \text{if $t \geq
1$}\end{cases}
$$
(because $\binom{2 n}{n} \leq 2^{2 n} = 4^n$ and $\binom{2 n +
t}{n} \leq 2^{2 n + t} = 2^t 4^n$). The lower bound of Corollary
\ref{c1} follows. $~~~~\blacksquare$ \vskip 1mm \noindent {\bf
Proof of Theorem \ref{t8}.} Theorem \ref{t8} is only a combination
of the results of Theorems \ref{t2} and \ref{t3} which we apply
for the arithmetic progression ${(u_{\ell})}_{\ell \in
\mathbb{N}}$ with general term $u_{\ell} = \ell + n$ (where $n \in
\mathbb{N}$ is fixed). $~~~~\blacksquare$ \vskip 1mm \noindent
{\bf Proof of Theorem \ref{t9}.} Let us prove the first assertion
of Theorem \ref{t9}. Giving $k$ a non-negative integer and $n$ a
positive integer, we easily show that for any non-negative integer
$j \leq k$, we have:
$$
n \binom{n + k}{k} \binom{k}{j} ~=~ (n + j) \binom{n + j - 1}{j}
\binom{n + k}{k - j} .
$$
It follows that the integer $\ppcm\!\left\{\!n \binom{n + k}{k}
\binom{k}{j} ~;~ j = 0 , \dots , k\!\right\} = n \binom{n + k}{k}
\ppcm\!\left\{\!\binom{k}{0} , \dots , \binom{k}{k}\!\right\}$ is
a multiple of each integer $n + j$ $(0 \leq j \leq k)$. Then it is
a multiple of $\ppcm\{n , n + 1 ,$\\$\dots , n + k\}$ as required.

Now, in order to prove the second assertion of Theorem \ref{t9},
we introduce the sequence of maps ${(g_k)}_{k \in \mathbb N}$ of
${\mathbb N}^*$ into ${\mathbb N}^*$ which is defined by:
$$g_k(n) := \frac{n (n + 1) \dots (n + k)}{\m\{n , n + 1 , \dots , n + k\}} ~~~~~~ (\forall k \in \mathbb N ,
\forall n \in {\mathbb N}^*) .$$ Let us show that ${(g_k)}_k$
satisfies the induction relation:
\begin{equation}
g_k(n) = \d\{k! , (n + k) g_{k - 1}(n)\} ~~~~~~ (\forall (k , n)
\in {{\mathbb N}^*}^2) . \label{1}
\end{equation}
For any pair of positive integers $(k , n)$, we have:
\begin{eqnarray*}
g_k(n) & := & \frac{n (n + 1) \dots (n + k)}{\m\{n , n + 1 , \dots , n + k\}}~\vspace{0.1cm} \\
& \!\!\!\!\!\!\!\!\!\!\!\!\!\!\!\!\!\!\!\!\!\!\!\!\!\!\!\!\!\!= &
\!\!\!\!\!\!\!\!\!\!\!\!\!\!\!\!\!
\frac{n (n + 1) \dots (n + k)}{\m\left\{\m\{n , n + 1 , \dots , n + k - 1\} , n + k\right\}}~\vspace{0.1cm} \\
& \!\!\!\!\!\!\!\!\!\!\!\!\!\!\!\!\!\!\!\!\!\!\!\!\!\!\!\!\!\!= &
\!\!\!\!\!\!\!\!\!\!\!\!\!\!\!\!\!\frac{n (n + 1) \dots (n +
k)}{\frac{\m\{n , n + 1 , \dots , n + k - 1\} .
(n + k)}{\d\left\{\m\{n , n + 1 , \dots , n + k - 1\} , n + k\right\}}}~\vspace{0.1cm} \\
& \!\!\!\!\!\!\!\!\!\!\!\!\!\!\!\!\!\!\!\!\!\!\!\!\!\!\!\!\!\!= &
\!\!\!\!\!\!\!\!\!\!\!\!\!\!\!\!\!\frac{n (n + 1) \dots (n + k -
1)}{\m\{n , n + 1 , \dots , n + k - 1\}} ~\!
\d\left\{\m\{n , n + 1 , \dots , n + k - 1\} , n + k\right\}~\vspace{0.1cm} \\
& \!\!\!\!\!\!\!\!\!\!\!\!\!\!\!\!\!\!\!\!\!\!\!\!\!\!\!\!\!\!= &
\!\!\!\!\!\!\!\!\!\!\!\!\!\!\!\!\!\d\left\{n (n + 1) \dots (n + k
- 1) , (n + k) g_{k - 1}(n)\right\} . \label{2}
\end{eqnarray*}
Then, the relation (\ref{1}) follows by remarking that the product
$n (n + 1) \dots (n + k - 1)$ is a multiple of $k!$ (because
$\frac{n (n + 1) \dots (n + k - 1)}{k!} = \binom{n + k - 1}{k} \in
\mathbb{N}$) and that $g_k(n)$ divides $k!$ (according to Theorem
\ref{t8}).

Now, giving a non-negative integer $k$, by reiterating the
relation (\ref{1}) several times, we obtain:
\begin{eqnarray*}
g_k(n) & = & \d\{k! , (n + k) g_{k - 1}(n)\} \\
& = & \d\{k! , (n + k) (k - 1)! , (n + k) (n + k - 1) g_{k - 2}(n)\} \\
& \vdots & \\
& = & \d\{k! , (n + k) . (k - 1)! , (n + k) (n + k - 1) . (k - 2)!
, \dots , \\
& ~ & ~~~~~~~~(n + k) (n + k - 1) \dots (n + k - \ell) g_{k - \ell
- 1}(n)\}
\end{eqnarray*}
for any positive integer $n$ and any non-negative integer $\ell
\leq k - 1$. In particular, for $\ell = k - 1$, since $g_0 \equiv
1$, we have for any positive integer $n$:
\begin{equation}
\begin{split}
g_k(n) &= \d\{k! , (n + k) . (k - 1)! , (n + k) (n + k - 1) . (k -
2)! , \dots ,\\
&~~~~~~~~~~~~(n + k) (n + k - 1) \dots (n + 1) . 0!\} \label{3}
\end{split}
\end{equation}
Now, if $n$ is a given positive integer satisfying the congruence
$n + k + 1 \equiv 0 ~\mod (k!)$, we have:
$$n + k \equiv - 1 ~\mod (k!) ~,~ (n + k) (n + k - 1) \equiv (- 1)^2 2! ~\mod (k!) ~,~ \dots ~,~$$
$$ (n + k) (n + k - 1) \dots (n + 1) \equiv (- 1)^k k! ~\mod(k!) ;$$
consequently, the relation (\ref{3}) gives:
$$g_k(n) = \d\left\{k! , 1! (k - 1)! , 2! (k - 2)! , \dots , k! 0!\right\} .$$
Hence:
\begin{eqnarray*}
\frac{k!}{g_k(n)} & = & \frac{k!}{\d\left\{0! k! , 1! (k - 1)! , \dots , k! 0!\right\}} \\
& = & \m\left\{\frac{k!}{0! k!} , \frac{k!}{1! (k - 1)!} , \dots , \frac{k!}{k! 0!}\right\} \\
& = & \m\left\{\binom{k}{0} , \binom{k}{1} , \dots ,
\binom{k}{k}\right\} .
\end{eqnarray*}
But on the other hand, according to the definition of $g_k(n)$, we
have:
$$\frac{k!}{g_k(n)} = \frac{\m\{n , n + 1 , \dots , n + k\}}{n \binom{n + k}{k}} .$$
We thus conclude that:
$$\m\{n , n + 1 , \dots , n + k\} = n \binom{n + k}{k} \m\left\{\binom{k}{0} , \binom{k}{1} ,
\dots , \binom{k}{k}\right\}$$ which gives the second assertion of
Theorem \ref{t9} and completes this proof.$~~~~\blacksquare$~\vspace{0.1cm}\\
{\bf Open Problem.} By using the relation (\ref{1}), we can easily
show (by induction on $k$) that for any non-negative integer $k$,
the map $g_k$ which we have introduced above is periodic of period
$k!$. In other words, the map $g_k$ $(k \in \mathbb{N})$ is
defined modulo $k!$. Then, for $k$ fixed in $\mathbb{N}$, it is
sufficient to calculate $g_k(n)$ for the $k!$ first values of $n$
$(n = 1 , \dots , k!)$ to have all the values of $g_k$.
Consequently, the relation (\ref{1}) is a practical mean which
permits to determinate step by step all the values of the maps
$g_k$. By proceeding in this way, we obtain: $g_0(n) \equiv g_1(n)
\equiv 1$ (obviously),
$$g_2(n) = \begin{cases}
1 & \text{if $n$ is odd} \\
2 & \text{if $n$ is even}
\end{cases} ~~,~~ g_3(n) = \begin{cases}
6 & \text{if} ~n \equiv 0 ~\mod(3) \\
2 & \text{otherwise}
\end{cases} ~~,~~ \dots\text{etc.}$$
This calculation point out that the smallest period of the map
$g_3$ is equal to $3 (\neq 3!)$. This lead us to ask the following
interesting open question:
\begin{center}
``Giving $k$ a non-negative integer, where is the smallest period
for the map $g_k$?''.
\end{center}

\end{document}